\documentclass[a4paper,12pt,reqno]{amsart}
\usepackage{amssymb}
\usepackage[dvips]{graphicx}
\usepackage{hyperref}

\input epsf.tex
\newdimen\xsize
\newdimen\oldbaselineskip
\newdimen\oldlineskiplimit
\def\restorelineskip{\baselineskip=\oldbaselineskip%
\lineskiplimit=\oldlineskiplimit}
\def\putm[#1][#2]#3{
\hbox{\vbox to 0pt{\parindent=0pt%
\vskip#2\xsize\hbox to0pt{\hskip#1\xsize $#3$\hss}\vss}}}%
\long\def\Line#1{\hbox to \hsize{#1}}
\def\putt[#1][#2]#3{
\vbox to 0pt{\noindent\hskip#1\xsize\lower#2\xsize%
\vtop{\restorelineskip#3}\vss}}

\makeatletter
\def\xbig[#1]#2{{\hbox{$\m@th\left#2\vbox to#1\xsize{}%
\right.\n@space$}}}
\makeatother
\def\xlar[#1]#2{%
\smash{\mathop{ \hbox to #1\xsize{\leftarrowfill}}\limits^{#2}}}
\def\xrar[#1]#2{%
\smash{\mathop{ \hbox to #1\xsize{\rightarrowfill}}\limits^{#2}}}
\def\xline[#1]{\hbox to #1\xsize{\leaders\hrule\hfill}}

\DeclareFontFamily{U}{rsf}{\skewchar\font'177}%
\DeclareFontShape{U}{rsf}{m}{n}{<-6>rsfs5<6-8>rsfs7<8->rsfs10}{}%
\DeclareFontShape{U}{rsf}{b}{n}{<-6>rsfs5<6-8>rsfs7<8->rsfs10}{}%
\DeclareMathAlphabet\RSFS{U}{rsf}{m}{n}
\SetMathAlphabet\RSFS{bold}{U}{rsf}{b}{n}
\def\sf#1{{\mathsf{#1}}}

\def\slsf{\slshape \sffamily }

\def\msmall#1{\mathchoice{\hbox{\small$\displaystyle {#1}$}}{#1}{#1}{#1}}

\let\xrar=\xrightarrow

\def\bb{{\mathbb B}}

\def\cc{{\mathbb C}}

\def\sss{{\mathbb S}}
\def\sph{{\mathbb S}}

\def\pp{{\mathbb P}}

\def\im{\sf{Im}\,}

\def\lim{\mathop{\sf{lim}}}

\def\w{{\mathrm{w}}}

\def\<{\langle}\let\la=\<
\def\>{\rangle}\let\ra=\>

\def\ddef{\mathrel{{=}\raise0.3pt\hbox{:}}}
\def\deff{\mathrel{\raise0.3pt\hbox{\rm:}{=}}}

\def\fraction#1/#2{\mathchoice{{\msmall{ #1\over#2}}}%
{{ #1\over #2 }}{{#1/#2}}{{#1/#2}}}

\def\le{\leqslant}

\def\longpoints{\leaders\hbox to 0.5em{\hss.\hss}\hfill \hskip0pt}
\def\stateskip{\smallskip}
\def\state#1. {\stateskip\noindent{\bf#1. }} 
\def\statep#1. {\stateskip\noindent{\bf#1 }} 
\def\proof{\state Proof. }

\def\Chi{\raise 2pt\hbox{$\chi$}}

\def\ie{\hskip1pt plus1pt{\sl i.e.\/,\ \hskip1pt plus1pt}}

\def\sli{{\sl i)} } 
\def\slii{{\sl i$\!$i)} } 
\def\sliii{{\sl i$\!$i$\!$i)} }

\def\Chi{\raise 2pt\hbox{$\chi$}}

\let\phI=\phi\let\phi=\varphi\let\varphi=\phI
\def\1{{1\mkern-5mu{\rom l}}}

\def\ge{\geqslant}

\def\fraction#1/#2{\mathchoice{{\msmall{ #1\over#2}}}%
{{ #1\over #2 }}{{#1/#2}}{{#1/#2}}}

\def\le{\leqslant}

\def\qed{\ \ \hfill\hbox to .1pt{}\hfill\hbox to .1pt{}\hfill $\square$\par}

\def\comment#1\endcomment{}

 \belowdisplayskip=\abovedisplayskip

\def\lineeqqno(#1){\hfill\llap{\vbox to 10pt%
{\vss\begin{align} \eqqno(#1)\end{align}\vss}}\vskip1pt}

\advance\headheight 1.2pt

\def\ShowwLLabel#1{}

\def\thechpt{\Roman{chpt}}

\def\newchapt[#1]#2{\newpage%
\refstepcounter{chpt}\setcounter{subsection}{0}%
\setcounter{thm}{0}\setcounter{defi}{0}%
\setcounter{rema}{0}\setcounter{exrc}{0}%
\renewcommand{\thesubsection}{\thechpt.\arabic{subsection}}%
\section*{\begin{center}\huge \bf Chapter \thechpt\\
#2 \end{center}}\label{#1}%
\ \smallskip%
\markboth{Chapter \thechpt}{#2}%
}
\def\newsect[#1]#2{\refstepcounter{section}\setcounter{equation}{0}%
\renewcommand{\thesubsection}{\arabic{section}.\arabic{subsection}}%
\section*{\arabic{section}.
#2}\vspace{-20pt}\label{#1}\vspace{20pt}%
\markboth{Section \arabic{section}}{#2}}

\def\newlect[#1]#2{\refstepcounter{section}%
\renewcommand{\thesubsection}{\arabic{section}.\arabic{subsection}}%
\section*{Lecture \arabic{section}\\
#2}\label{#1}%
\markboth{Lecture \arabic{section}}{#2}}

\def\newprg[#1]#2{\refstepcounter{subsection}%
\subsection*{{\thesubsection.\ #2}} \label{#1}%
}

\def\newappx[#1]#2{%
\refstepcounter{appx}\setcounter{section}{0}%
\renewcommand{\thesubsection}{A\arabic{appx}.\arabic{subsection}}%
\section*{Appendix \arabic{appx}\\ #2}
\label{#1}%
\markboth{Appendix A\arabic{appx}}{#2}
}

\newtheorem{thm}{Theorem}[section]
   \def\newthm#1{\begin{thm}\label{#1}}

\newtheorem{nnthm}{Theorem}
   \def\newthm#1{\begin{nnthm}\label{#1}}

\newtheorem{lem}{Lemma}[section]
   \def\newlemma#1{\begin{lem} \label{#1}}

\newtheorem{prop}{Proposition}[section]
   \def\newprop#1{\begin{prop}\label{#1}}

\newtheorem{nnprop}{Proposition}
   \def\newprop#1{\begin{nnprop}\label{#1}}

\newtheorem{corol}{Corollary}[section]
   \def\newcorol#1{\begin{corol} \label{#1}}

\newtheorem{nncorol}{Corollary}
   \def\newcorol#1{\begin{nncorol} \label{#1}}

\newtheorem{defi}{Definition}[section]
   \def\newdefi#1{\begin{defi} \label{#1}\rm }

\newtheorem{exmp}{Example}[section]
   \def\newexmp#1{\begin{exmp} \label{#1}\rm }

\newtheorem{nnexmp}{Example}
   \def\newexmp#1{\begin{nnexmp} \label{#1}\rm }

\newtheorem{exrc}{Exercise}
   \def\newexrc#1{\begin{exrc} \label{#1}\rm }

\newtheorem{rema}{Remark}[section]
   \def\newrema#1{\begin{rema} \label{#1}\rm }

\newtheorem{nnrema}{Remark}
   \def\newrema#1{\begin{nnrema} \label{#1}\rm }

\def\eqqno(#1){\label{(#1)}}
\def\eqqref(#1){(\ref{(#1)})}

\title{An example in concern with holomorphicity\\ 
of meromorphic mappings along\\ real hypersurfaces}
\author{S. Ivashkovich, F. Meylan}
\date{\today}

\address{
Universit\'e de Lille-1, UFR de Math\'ematiques, 59655 Villeneuve
d'Ascq, France} \email{ivachkov@math.univ-lille1.fr}
\address{IAPMM Nat. Acad. Sci. Ukraine
Lviv, Naukova 3b, 79601 Ukraine}

\address{Universit\'e de Fribourg, D\'epartement de Math\'ematiques,
1700 Fribourg, Suisse} \email{francine.meylan@unifr.ch}

\subjclass[2010]{Primary - 32H04, Secondary - 32V99} \keywords{
Meromorphic mapping, real hypersurface.}

\begin{document}
\begin{abstract}
We construct an example of a rational mapping $F$ from $\cc^2$ to $\pp^2$,
which has indeterminacies on the unit sphere $\sph^3\subset \cc^2$, but
such that $F|_{\sss^3}$ is continuous, image $K\deff F(\sph^3)$ is contained 
in the affine part $\cc^2$ of $\pp^2$ and $K$ doesn't contain any germ of a 
non-constant complex curve.
\end{abstract}

\maketitle

\setcounter{tocdepth}{1}
\tableofcontents

\newsect[sect.INT]{Introduction.}

Let $F:U\to \pp^N$ be  a meromorphic mapping  from a domain $U\subset \cc^n$ 
to the complex projective space. Here $n\ge 2$ and $N\ge 2$ will be supposed 
everywhere. Denote by $I_F$ the set of points of indeterminacy of $F$, \ie 
$z_0\in I_F$ if and only if $F$ is not holomorphic in any neighborhood of 
$z_0$. As it is well known $I_F$ is an analytic subset of $U$ of codimension 
at least two. Recall that the {\slsf full image} by $F$ of a point $z_0$, 
denoted as $F[z_0]$, is the set of all cluster points of $F$ at $z_0$, \ie 
\begin{equation}
\eqqno(full-im)
F[z_0] = \bigl\{x\in \pp^N: \exists z_k\to z_0, z_k\notin I_F \text{ such that } 
F(z_k)\to x\bigr\}.
\end{equation}
$F$ is holomorphic in a neighborhood of $z_0$ if and only if $F[z_0]$ is a
singleton. Likewise one can define the full image of $z_0$ by $F$ along a 
closed subset $M$ accumulating to $z_0$, ex. $M$ can be a complex curve or,
a real hypersurface containing $z_0$:
\begin{equation}
F_M[z_0] = \bigl\{x\in \pp^N: \exists m_k\in M\setminus I_F, m_k\to z_0 
\text{ such that } F(m_k)\to x\bigr\}.
\end{equation}
If $M$ is a complex curve then $F_M[z_0]$ is always a singleton regardless 
of whether $z_0$ is a regular or an indeterminacy point of $F$. As we shall 
see in our example the same can happen if $M$ occurs to be a real hypersurface
with $z_0$ being an indeterminacy point of $F$. The proper image or transform
of a closed subset $M\subset U$ under $F$ is defined now to be the union of 
full images along $M$ of its points. In other words 
\begin{equation}
\eqqno(prop-im1)
F_M[M] = \bigcup_{m\in M}F_M[m].
\end{equation}
If $M$ is compact, for example $M=\sph^3\subset \cc^2 = U$, then 

\begin{equation}
\eqqno(prop-im2)
F_M[M] = \overline{F(M\setminus I_F)}.
\end{equation}

We denote as $w_1,w_2$ the standard coordinates in $\cc^2$. $\sph^3=\{(w_1,w_2):
|w_1|^2 + |w_2|^2=1\}$ stands for the unit sphere in $\cc^2$ and $\bb^2$ for the 
{\slsf open} unit ball. Our goal in this note is to construct the following example:

\begin{nnexmp}
There exists a rational mapping $F:\cc^2\to \pp^2$ such that:

\smallskip\sli Its indeterminacy set is $I_F=\{p_{\pm}\deff (0,\pm 1)\}$,
its divisor of poles is $P_F = \{w_2 = \pm 1\}$, 

\quad and $F_{\sph^3}[p_{\pm}]$ are singletons, \ie $F|_{\sss^3}$ is continuous.

\smallskip\slii The proper image $K\deff F_{\sph^3}[\sph^3]$ of the unit 
sphere $\sph^3$ under $F$ doesn't contain any 

\quad germ of a non-constant complex curve.

\smallskip\sliii In addition $K$  is contained in an appropriate affine part $W\cong\cc^2$ of 
$\pp^2$.
\end{nnexmp}

Let us explain the interest in such an example. It was proved in  Lemma 6.6, p. 128 
of \cite{P3} and later, independently in \cite{C}, that for a germ $M$ of a real
analytic hypersurface in $\cc^n$ and a germ $F$ of a meromorphic mapping near $M$ 
such that $F_M[M]\subset \sph^{2N-1}$, $2\le n\le N$, one has that $F$ is holomorphic 
near $M$. It is of considerable interest and importance to replace in this result
the unit sphere $\sph^{2N-1}\subset \cc^N$ by an arbitrary compact real
analytic set $K\subset \cc^N$ which doesn't contain non-constant germs of
complex curves, for instance  $K$  smooth and  strictly pseudoconvex.  This problem 
has been under discussion among the experts, and our example shows that the conjecture 
mentioned above can not hold if one does not  assume that $K$ is smooth. See the 
relevant Remarks \ref{rasul} and \ref{eq-dim} on these issues at the end of this note.

\smallskip\noindent{\slsf Acknowledgements.} We are grateful to Rasul Shafikov and Alexander
Sukhov for the stimulating discussions, in particular for the Remark \ref{rasul}
about our example and for giving us the possibility to consult the thesis of 
S. Pinchuk \cite{P3}.

\newsect[PR]{Construction of the Example}

\smallskip\noindent{\slsf Step 1. \it Construction of the mapping.}
As it is well known the unit sphere $\sph^3\subset\cc^2\subset\pp^2$ after an appropriate
projective transformation of $\pp^2$, see \eqqref(transf) below, in an appropriate affine chart 
$U_0=\cc^2$ can be represented in the form
\begin{equation} 
\eqqno(eqm0)
M \deff \{\im z_2 =  |z_1|^2\}.
\end{equation}

\smallskip As a rational map $F$ in question we shall take 
\begin{equation}
\eqqno(map-f)
F :(z_1,z_2) \to \left(\frac{z_1^3}{z_2},z_2\right).
\end{equation}
For the moment it is better to consider this map as a mapping from $U_0 = \cc^2$ 
to $\pp^1\times \cc$. Therefore $F$ is holomorphic on $M\setminus \{0\}$. To see
that $F$ is continuous on the whole of $M$ one should remark that  $F_M[0] = 0$, 
\ie the full image of $0$ along $M$ is a singleton. This is nearly obvious because
along $M$ one has $|z_2|\ge |z_1|^2$ and therefore along $M$ the fist function
in \eqqref(map-f) is $O(|z_1|)$.

\smallskip Set $l\deff M\cap \{z_1=0\}$, this is a real line. One easily checks 
that $F$ is identity when restricted to $l$ and that on $M\setminus l$ it is an
injection. The latter can be seen from computing the Jacobian of $F$. 
$F|_{M\setminus l}:M\setminus l\to F(M\setminus l)$ is therefore a three to one 
cover. This means that $F(M\setminus l)$ is a strictly pseudconvex hypersurface.
Therefore it cannot contain any germ of a complex curve. As a result a germ of
a complex curve in $F(M)$ should be contained in $F(l)=l$, which is impossible.
I.e., $F(M)$ doesn't contain germs of comlex curves.

\medskip\noindent{\slsf Step 2. \it Globalization.} It is easy 
to check that the projective transformation 
\begin{equation}
\eqqno(transf)
R = 
\begin{cases}
z_1 = \frac{w_1}{1-w_2}\cr
z_2 = i\frac{1+w_2}{1-w_2}
\end{cases}
\end{equation}
sends the unit sphere $\sph^3 = \{|w_1|^2+|w_2|^2 = 1\}$ to our model $M$. Writing 
$F$ as a mapping to $\pp^2$, \ie as 
\begin{equation}
\eqqno(map-f-p2)
F: (z_1,z_2) \to [z_1^3:z_2^2:z_2],
\end{equation}
and composing it with $R$ we get the expression of our map in coordinates $w_1,w_2$:
\begin{equation}
\eqqno(map-f-r)
F : (w_1,w_2) \to [w_1^3:-(1+w_2)^2(1-w_2):i(1+w_2)(1-w_2)^2].
\end{equation}
We refer to \eqqref(map-f-r) as to the form of the mapping $F$ in coordinates
$(w_1,w_2)$. This is the mapping of our example. Now let us prove the 
assertions (\sli - (\sliii as they are stated in the Introduction.

\smallskip\noindent (\sli It is easy to see that it 
has two points of indeterminacy $p_{\pm} = (0,\pm 1)$ and that it behaves
at this points similarly. The behavior at $p_-$ was already studied,
the behavior at $p_+$ is the same. Divisor of poles is $\{w_2=\pm 1\}$ as 
required. 

\smallskip\noindent (\slii $K:=F(\sss^2)$ cannot contain any germ of a non-constant 
complex curve. Indeed, let $C\subset K$ be such a germ. Then $C$ should 
intersect an affine neighborhood either of $p_+$ or, of $p_-$, and this was proved to be
not the case.

\smallskip\noindent (\sliii To prove that $K$ is situated in an appropriate affine part
of $\pp^2$ consider in homogeneous coordinates $[\zeta_0:\zeta_1:\zeta_2]$ 
the line 
\begin{equation}
\eqqno(line-l)
L\deff \{\zeta_2 - i\zeta_1 =0\} .
\end{equation}
Putting $F(w)$ to this equation we get 
\begin{equation}
\eqqno(f-to-line)
\zeta_2 - i\zeta_1 = i(1+w_2)(1-w_2)^2 + i(1+w_2)^2(1-w_2) = 2i(1-w_2^2).
\end{equation}
The latter is not zero on $\sph^3$ if and only if $\w_2\not=\pm 1$, \ie if $w$ is 
not $p_{\pm}$. But at $p_{-}$ along $\sph^3$ our map takes value $[0:0:1]\not\in
L$ and at $p_+$ the value $[0:1:0]\not\in L$. Therefore $F_{\sph^3}[\sph^3] \subset 
\pp^2\setminus L =:W\cong\cc^2$. Now we can can make a coordinate change 
\begin{equation}
\begin{cases}
\xi_0 = \zeta_0\cr
\xi_1 = \zeta_1\cr
\xi_2 = \zeta_2 - i\zeta_1
\end{cases}
\end{equation}
in $\pp^2$ and write our map in the form
\begin{equation}
\eqqno(map-f-final)
F : (w_1,w_2) \to [w_1^3:-(1+w_2)(1-w_2^2):2i(1-w_2^2)],
\end{equation}
or, in the affine charts 
\begin{equation}
\eqqno(map-f-affine)
F : (w_1,w_2) \to \left(-\frac{i}{2}\frac{w_1^3}{1-w_2^2}, \frac{i}{2}(1+w_2)\right).
\end{equation}
This is the final form of $F$ as a mapping from $\cc^2\subset\pp^2$ to $\cc^2\subset \pp^2$,
\ie $K\subset \cc^2$ id $F$ is given in the form \eqqref(map-f-affine).

\smallskip Let us make few additional remarks about our example.

\begin{rema} \rm 
\label{rasul}

\noindent{\bf 1.} Remark that the image $K=F(\sss^2)$ in our example is not smooth.
Indeed it is given by the equation 
\begin{equation}
\eqqno(K-eq)
(\im z_2)^3 = |z_1z_2|^2.
\end{equation}
It is easy to see that the set of singular points of \eqqref(K-eq) is exactly our
line $l=F(M)\cap \{z_1=0\}$ considered above. Outside of $l$ our $K$ is a spherical
hypersurface.

\smallskip\noindent{\bf 2.} Non-smoothness of $K$ can be overcomed on the price of
incresing the dimension of the image. Our example written as a map $\cc^2\to \cc^3$, \ie as 
\begin{equation}
F:(z_1,z_2) \to \left(\frac{z_1^3}{z_2},z_2,0\right),
\end{equation}
takes its values in the smooth algebraic hypersurface 
\begin{equation}
K = \{y_3 + y_2^3 - |z_1z_2|^2 = 0\}.
\end{equation}
But this hypersurface contains the complex curve $C = \{z_2 = z_3 =0\}$.
\end{rema}

\begin{rema} \rm
\label{eq-dim}
In equal dimensions with strictly pseudoconvex both $M$ and $K$ such example doesn't
exist. Namely the following is true.

\begin{prop}
\label{eq-dilm-case}
Let $(M,0)$ be a germ of strictly pseudoconvex hypersurface in $\cc^n$, $n\ge 2$, 
and let $K$ be a compact strictly pseudoconvex hypersurface in $\cc^n$. Let  
$F:V\to \cc^n$ be a meromorphic mapping from some neighborhood of the origin 
$V\subset \cc^n$ to $\cc^n$ such that $F(M)\subset K$. Then $F$ is holomorphic
in a neighborhood of the origin.
\end{prop}
\proof If $F$ is constant then there is nothing to prove. If, not then it is 
locally biholomorphic at every point of $M\setminus I_F$. Now one has two cases.

\smallskip\noindent{\slsf Case 1. \it $M$ is spherical.} Then $K$ is spherical
as well. Set $p\deff F(0)$ and let $j:(K,p)\to (\sss^{2n-1},0)$ be a germ of 
a biholomorphic mapping. By the result of  \cite{P1} germ $j$ extends locally
biholomorphically along any path in $K$. Take some point $a\in M\setminus I_F$
and set $b\deff F(a)$. Take a path $\gamma$ in $M$ from $0$ and $a$. Set $\beta
\deff F(\gamma)$, it is a path from $p$ to $b$ in $K$. Composition $j\circ F$
extends along $\gamma^{-1}$ locally biholomorphically form $a$ to $0$ by the 
mentioned result of \cite{P1}. 
Since $j$ extends along $\beta$ and therefore $j^{-1}$ extends along $j(\beta^{-1})$
we can take a composition $j^{-1}\circ (j\circ F)$ and extend it from $a$ to $0$
to obtain the holomorphicity of $F$ at the origin.

\smallskip\noindent{\slsf Case 2. \it $M$ is not spherical.} Then $K$ is not
spherical as well. We can then apply another Theorem of Pinchuk, see \cite{P2}.
Take any Moser chain through $0$ on $M$ which is not contained entirely in 
$M\cap I_F$, say $\gamma:[0,1]\to M$ such that $\gamma (0)\not\in I_F$ and
$\gamma (1)=0$. Extend $F$ along $\gamma$ by \cite{P2} to conclude that
$F$ is holomorphic at zero.

\smallskip\qed
\end{rema}

\ifx\undefined\bysame
\newcommand{\bysame}{\leavevmode\hbox to3em{\hrulefill}\,}
\fi

\def\entry#1#2#3#4\par{\bibitem[#1]{#1}
{\textsc{#2 }}{\sl{#3} }#4\par\vskip2pt}

\end{document}